\documentclass[11pt]{article}
 
\begin{document} 

\textwidth  7in
\textheight 8.5in

\title{Trajectory and Attractor Convergence for 
a Nonlocal Kuramoto-Sivashinsky Equation}
\author{  }
\date{October 25, 1997 }
\maketitle

\begin{center}

Jinqiao Duan$^{1}$, Vincent J. Ervin$^{1}$ and Hongjun Gao$^{2}$  \\

{\em 
1. Department of Mathematical Sciences, Clemson University, \\
   Clemson, South Carolina 29634, USA. \\
  
2. Laboratory of Computational Physics, \\
   Institute of Applied Physics and Computational Mathematics, \\
   P. O. Box 8009, Beijing 100088, China. \\}
 
\end{center}

\newcommand{\e}{\epsilon}
\renewcommand{\a}{\alpha}
\renewcommand{\b}{\beta}
\renewcommand{\aa}{\mbox{$\cal A$}}
\renewcommand{\L}{\mbox{$\cal L$}}
 
\renewcommand{\theequation}{\thesection.\arabic{equation}}
\newtheorem{theorem}{Theorem}[section]
\catcode`\@=11
\@addtoreset{equation}{section}

\begin{abstract}

The nonlocal Kuramoto-Sivashinsky equation arises 
in the modeling of the flow
of a thin film of viscous liquid falling down an inclined plane,
subject to an applied electric field. 
In this paper, the authors show that, as
the coefficient of the nonlocal integral term goes to zero,
the solution trajectories and the maximal attractor
of the nonlocal Kuramoto-Sivashinsky equation converge  to
those of the usual Kuramoto-Sivashinsky equation.



\bigskip
\bigskip

{\bf To Appear:} Communications on Applied Nonlinear Analysis

		 No. 4, Vol. 5, 1998.
\end{abstract}
 
\newpage

	      \section{Introduction}

In this paper we investigate some dynamical issues of a  
Kuramoto-Sivashinsky equation with a nonlocal spatial integral term
\begin{equation}
u_t + u_{xxxx} +u_{xx} + u u_x +\alpha H(u_{xxx}) =0,
\label{ks}
\end{equation}
where $\alpha$ is a non-negative constant coefficient, and
\begin{equation}
 H(f) =\frac{1}{\pi} \int_{-\infty}^{\infty} \frac{f(\xi)}{x-\xi}d\xi,  
\label{hil1}
\end{equation}
is the Hilbert transform and the integral is understood in the sense 
of the Cauchy Principle Value.
 
The equation under consideration arises in the modeling of the flow
of a thin film of viscous liquid falling down an inclined plane,
subject to an applied electric field \cite{nonlocal_ks}. The
application of a uniform
electric field at infinity, perpendicular to the inclined plane,
is to destabilize the liquid films on the surface of the plane.
In an industrial setting it is hoped that this destabilization
will lead to an enhancement of heat transfer.
In this engineering context,  
$\alpha W_{e} / \sqrt{W |\cot \beta - \frac{4}{5} R_{e} |}$,
where $W_{e}$ is the electrical Weber number, $W$ the Weber number,
$R_{e}$ the Reynolds number and $\beta$ the
angle between the plane and the horizontal.
When the electrical field is absent, i.e. $W_{e} = 0$, 
the Hilbert transform term is
gone and we obtain the usual (local) Kuramoto-Sivashinsky equation 
\begin{equation}
u_t + u_{xxxx} +u_{xx} + u u_x =0.
\label{ks1}
\end{equation}

Throughout this paper we restrict our attention to the case of 
$u(x,t)$ periodic in $x$ with period $l>0$.
We denote the interval $I := (-l, l)$. 
Then, (\ref{hil1}) is replaced by 
(see \cite{Bona}) the ``periodic"
Hilbert transform
\[ H(f) = -\frac{1}{2l}\int_{I} cot\frac{\pi (x-\xi)}{2l}\; f(\xi)d\xi .\]

Observe that $u = C$, a constant, satisfies (\ref{ks}). Moreover, for
\begin{equation}
     u(x , 0) \, = u_{0}(x) \, \mbox{ for } x \in I \, , \mbox{ and }
       \bar{u}_{0} \, := \, \frac{1}{2 l} \int_{I} u_{0}(x) \, dx \, ,
\label{skb1}
\end{equation}
integrating (\ref{ks}) over the interval I yields
\[ 
     \frac{d}{d t} \int_{I} u(x , t) \, dx \, = \, 0 \, , \, \mbox{ i.e. } 
     \frac{1}{2 l} \int_{I} u(x , t) \, dx \, = \, \bar{u}_{0} \, .
\]
Thus, we have that the mean of the solution is conserved with respect to 
time. This may be interpreted in the sense that the ``dynamics'' of $u$
satisfying (\ref{ks}) are centered around the mean value of the initial data. 
Therefore, without loss of generality, we will work in a periodic functional
space with zero mean. That is, we will  assume the following condition
\begin{equation}
 u(x, t) \mbox{ periodic on } I, \mbox{ with } u(x, 0) = u_{0}(x) \, ,
\mbox{ and } \int_{I} u_{0}(x) \, dx \, = \, 0 \, .
\label{bc}
\end{equation}

In this paper, we will show that individual solution trajectories and 
the maximal attractor of the nonlocal Kuramoto-Sivashinsky equation approach
the solution trajectories and the maximal attractor of the
usual local Kuramoto-Sivashinsky equation, as the coefficient for
the nonlocal term goes to zero.

		\section{Preliminaries}

We denote by $L^2_{per}(I)$, $H^k_{per}(I),\; k=1,2,\cdots$, 
the usual Sobolev
spaces of periodic functions on $I$. Let $\dot{L}^2_{per}(I)$, $\dot{H}^k_{per}(I)$
denote the spaces of 
functions $g$ in  $L^2_{per}(I)$, $H^k_{per}(I)$, respectively, with
mean zero, i.e. $\bar{g} := \frac{1}{2 l} \, \int _{I} g(x)dx \, = \, 0$.
In the following, $||\cdot ||$ denotes the usual $L^2_{per}(I)$ norm.
Due to the Poincare inequality, $\|D^k u\|$ is an equivalent norm in
$\dot{H}^k_{per}(I)$.
All integrals are with respect to $x\in I$, unless specified otherwise.
We recall the following two inequalities:

{\em Poincar\'{e} inequality} (\cite{Temam88}, p.49):
\begin{equation}
\int_{I} g^2 \, dx \leq (2l)^2 \int_{I} g_{x}^2 \, dx,
 \label{peieq}
\end{equation}
for $g \in \dot{H}^{1}_{per}(I)$.

{\em Agmon inequality} (\cite{Temam88}, p.50):
\begin{equation}
||g||_{\infty}^2 \leq 2 ||g|| \, ||g_x||,
 \label{anieq}
\end{equation}
for $g \in \dot{H}^1_{per}(I)$. This inequality also follows from
$g^2(x) = 2\int_{x_0}^x g g_x \,dx$, where $g(x_0) = 0$ (as
$g$ is continuous and has zero mean).

The Hilbert transform (\ref{hil1}) is a linear, invertible,
bounded operator
from $L^2$ to $L^2$, and from Sobolev space $H^k$ to  $H^k$.
Several noteworthy properties of the transform are 
(see \cite{Bona}, \cite{Alarcon} or \cite{Stein} ):
\begin{eqnarray}
D_x H &=& H (D_x)    \, ,  \label{HDex} \nonumber\\
H^{-1} &=& -H, \nonumber \\
\int v H(u) &=& -\int u H(v), \nonumber \\
\int H(u) H(v) &=& \int u v,  \nonumber \\
\int u H(u) &=& 0, \label{equuHu} \nonumber \\
\| H(u) \| &=& \| u \|,   \nonumber \\
(H(u))_x  & = & H(u_x).   \nonumber
\end{eqnarray}

These properties hold for the Hilbert transformation
on both the real line and periodic intervals (\cite{Bona}).
On the periodic interval $(-l , l)$, the Hilbert transformation has
a simple representation
\[
H(f)(x) = i \sum_{k\in Z} sgn(k)f_k e^{ik\pi x/l},
\]
for $f(x) = \sum_{k\in Z} f_k e^{ik\pi x/l}$, with $f_k$'s the Fourier
coefficients of $f$.

Duan and Ervin \cite{Duan_Ervin} (or, for similar
arguments, see \cite{DuanLT1}, \cite{DuanLT2} or \cite{Temam88}), 
have shown that the nonlocal dynamical system (\ref{ks}) and (\ref{bc})
has global solutions and a maximal attractor in $\dot{L}^2_{per}(I)$. 
Moreover the following time-uniform estimates hold
for the nonlocal Kuramoto-Sivashinsky equation:
\begin{eqnarray}
  \|u\|_{L^2_{per}} \leq \rho_0, \\
  \|u\|_{H^1_{per}} \leq \rho_1,  \label{rho1} \\
  \|u\|_{H^2_{per}} \leq \rho_2,  \label{rho2}
\end{eqnarray}
where $\rho_0, \rho_1, \rho_2$ are absolute positive constants
for $0 \leq \alpha \leq 1$. These  estimates
have also been proved for the usual Kuramoto-Sivashinsky equation   
(\cite{Collet}, \cite{Temam_Wang} and references therein).
 
The usual Kuramoto-Sivashinsky equation (\ref{ks1}) together
with (\ref{bc}) also has 
global solutions and a maximal attractor in $\dot{L}^2_{per}(I)$
as well as other dynamical properties; see, for example,
\cite{NST85}, \cite{NST89}, \cite{Ilynko}, \cite{Goodman},
\cite{Jolly} and \cite{Robinson}.

\section{Localization Limit}

In this section, we will show that individual solution trajectories and
the maximal attractor for the nonlocal  
Kuramoto-Sivashinsky equation (\ref{ks})
converge to the solution trajectories and the 
maximal attractor of the usual Kuramoto-Sivashinsky equation (\ref{ks1}),
as the coefficient $\alpha$ of the nonlocal term goes to zero. We note that
Duan et al. \cite{DuanLT2} have studied the this convergence for a
Kuramoto-Sivashinsky equation with an extra dispersive term $u_{xxx}$.
Ercolani et al. \cite{EMR} have also studied the impact of $u_{xxx}$
on the dynamics of the usual Kuramoto-Sivashinsky equation.

General convergence result for the maximal  attractors is available in
Hale et al. \cite{Hale}, Temam \cite{Temam88} and Hale-Lin-Raugel
\cite{Hale2}. Similar applications can also be found in
Hill and Sulin \cite{Hill}.

Let us denote by $T_{\alpha}(t)$ and $T(t)$ the solution operators
for the nonlocal Kuramoto-Sivashinsky equation (\ref{ks}),
and the usual Kuramoto-Sivashinsky equation (\ref{ks1}), respectively.
We need to verify a few conditions that will ensure   the
convergence of the solution trajectories and attractors. 
In particular, we should show that operators 
$T_{\alpha}(t)$ and $T(t)$ satisfy
\[
\|T_{\alpha}(t)u_0 - T(t)u_0\| \le \eta(\alpha,t,u_0),
\]
where $ \eta(\alpha,t,u_0) \rightarrow 0$ as $\alpha \rightarrow 0$,
and that the domain of attraction of $T_{\alpha}(t)$ is independent of $\alpha$.
For more details of these conditions see Theorem I.1.2 of Temam \cite{Temam88},
or Theorem 2.4 of Hale-Lin-Raugel \cite{Hale2}.

Assume that $u, v$ are the solutions of the nonlocal and usual
Kuramoto-Sivashinsky equations, respectively, 
with the same initial data $u_0(x)$,
i.e., $u=T_{\alpha}(t)u_0$, $v=T(t)u_0$. We denote also
$w=u-v$, and hence $\|w(0)\| \equiv 0$.

Using the equations (\ref{ks}) and (\ref{ks1}), we get
\begin{eqnarray}
\frac12 \frac{d}{dt} \|w\|^2 + \|w_{xx}\|^2 -\|w_x\|^2
&+& \int (uw_x+wv_x)wdx        \nonumber \\
&+&\alpha \int H(w_{xxx})wdx =0.
\label{w}
\end{eqnarray}
After integration by parts, and using the Cauchy-Schwarz inequality and
Agmon inequality (\ref{anieq}), we get
\begin{eqnarray}
\|w_x\|^2
& \leq &  \|w\|^2 +\frac14 \|w_{xx}\|^2,
\label{w1}  \\
\int (uw_x+wv_x)wdx
& = &  \int (-\frac12 u_x +v_x) w^2 dx   	\nonumber \\
& \leq & (\frac12 \|u_x\|^{\frac12} \|u_{xx}\|^{\frac12}
	+\|v_x\|^{\frac12} \|v_{xx}\|^{\frac12})\|w\|^2    \nonumber \\
& \leq & 2\rho_1^{\frac12} \rho_2^{\frac12}\|w\|^2,
\label{w2}   \\
\alpha \int H(w_{xxx})wdx = -\alpha \int H(w_{xx})w_x dx
& \leq &  \alpha^2 \|w_{x}\|^2
	+\frac14 \|w_{xx}\|^2           \nonumber \\
& \leq &   \alpha^2  \rho_1^2 + \frac14 \|w_{xx}\|^2.
\label{w3}
\end{eqnarray}
Here we have used the estimates in (\ref{rho1}) and (\ref{rho2}).
Substituting (\ref{w1}), (\ref{w2}) and (\ref{w3}) into (\ref{w}), we obtain,
\begin{eqnarray}
\frac{d}{dt} \|w\|^2 + \|w_{xx}\|^2 \leq  2 \|w\|^2
+ 4 \rho_1^{\frac12} \rho_2^{\frac12}\|w\|^2 + 2\alpha^2 \rho_1^2.
\label{w4}
\end{eqnarray}
Let
\[     b_1 =  2 + 4 \rho_1^{\frac12} \rho_2^{\frac12}, \:\: \mbox{\rm
and}
 \:\: b_2= 2 \rho_1^2.   \]
Thus we have
\begin{eqnarray}
\frac{d}{dt} \|w\|^2 \leq b_1 \|w\|^2 + b_2 \alpha^2,
\end{eqnarray}
or
\begin{eqnarray}
\|w(t)\|^2 \leq \|w(0)\|^2 e^{b_1 t} + \frac{b_2}{b_1}(e^{b_1 t}-1)\alpha^2,
\end{eqnarray}
Since $\|w(0)\|\equiv 0$, we conclude that $\|w(t)\| \rightarrow 0$ as
$\alpha \rightarrow 0$. This shows that the solution 
$u$ of the nonlocal Kuramoto-Sivashinsky equation (\ref{ks}) approaches
the solution $v$ of the usual Kuramoto-Sivashinsky equation (\ref{ks1}).

Moreover, solution operator $T_{\alpha}(t)$ has a uniform
domain of attraction, i.e., $\dot L^2_{per}$, which is independent of
$\alpha$.

Therefore, by using Theorem I.1.2 of Temam \cite{Temam88}
or Theorem 2.4 of Hale-Lin-Raugel \cite{Hale2},
we have the following theorem.

\begin{theorem}
When the coefficient $\alpha$ of the nonlocal integral term goes to zero,
 
(i) A solution trajectory
of the nonlocal Kuramoto-Sivashinsky equation (\ref{ks}) converges in $L^2$ to
a solution trajectory of the usual Kuramoto-Sivashinsky equation (\ref{ks1}),
as long as both trajectories start at the same initial point.

(ii) The maximal attractor of the the nonlocal Kuramoto-Sivashinsky equation  
(\ref{ks}) converges in $L^2$
to the maximal attractor of the usual Kuramoto-Sivashinsky equation
(\ref{ks1}).
\label{convergence}
\end{theorem}

	\section{Discussions}

      Nonlocal spatial integral terms appear in mathematical models in the
form of partial integro-differential equations for various physical
phenomena, for example, in beam theory \cite{Holmes},
in electric systems \cite{Radehaus1990}, 
chemical diffusion \cite{Graham},
heat conduction processes \cite{Lacey},
temperature evolution in an atmospheric system \cite{Murthy}, 
in ferromagnetic system \cite{Elmer},
and in chemical frontal development \cite{Wilder}.
 	
	In this paper we have studied the limiting behavior of a 
prototypical partial-integral differential equation arising in fluid
dynamics, as the coefficient of the nonlocal integral term goes to zero. 
This work is a part of our effort on the investigation of
impact of nonlocal spatial interactions on the dynamics of
(local) infinite dimensional dynamical systems.

\bigskip
 
{\bf Acknowledgement}. This work was supported by the U.S.
National Science Foundation grant DMS-9704345, and by the National 
Natural Science Foundation of China grant 19701023.

\end{document}